\documentclass[letterpaper, 10 pt, conference]{ieeeconf}  

\IEEEoverridecommandlockouts                              
\title{\LARGE \bf
Dual  Particle Output Feedback Control based on Lyapunov drifts for nonlinear systems
}

\usepackage[autostyle]{csquotes}
\usepackage{amsmath,verbatim,amssymb,amsfonts,amscd, graphicx}
\usepackage{mathabx}
\usepackage{breqn}
\usepackage{float}
\usepackage{graphics}
\usepackage{stmaryrd} 
\usepackage{textcomp}
\usepackage{setspace}
\usepackage{algorithm}
\usepackage{algpseudocode}

\author{Emilien Flayac, Karim Dahia, Bruno H\'eriss\'e, and Fr\'ed\'eric Jean}

\begin{document}

\maketitle
\thispagestyle{empty}
\pagestyle{empty}

\begin{abstract}
This paper presents a dual receding horizon output feedback controller for a general non linear stochastic system with imperfect information. The novelty of this controller is that stabilization is treated, inside the optimization problem, as a negative drift constraint on the control that is taken from the  theory of stability of Markov chains. The dual effect is then created by maximizing information over the stabilizing controls which makes the global algorithm easier to tune than our previous algorithm. We use a particle filter for state estimation to handle nonlinearities and multimodality. The performance of this method is demonstrated on the challenging problem of terrain aided navigation.  

\end{abstract}

\section{INTRODUCTION}

Stochastic model predictive control (SMPC) is a widespread technique to deal with control problems where the state is subject to stochastic disturbances, hard constraints on the input and potentially soft constraints on the state. Its principle is to apply a receding horizon strategy based on the resolution of a finite horizon stochastic optimal control problem. Although, when only partial and noisy information on the state is available through some observations, classic SMPC combined with a state estimator may lead to overcautious controls or even destabilizing ones. It is due to the fact that, in general, the control influences the observations and then state estimation in addition to guiding the system in a standard way. This property is known as the \textit{dual effect} property of the control. Stochastic optimal control problems with imperfect information are much harder than their counterpart in the perfect information case. The main reason is that to be optimal with partial information, one needs to anticipate the information that will be available. Consequently, the optimal law possesses the property of dual effect. As it is computationally intractable, suboptimal outputfeedback control laws are computed instead, with the idea to keep the dual effect property. Therefore, they are called dual controllers. 

Usually, dual controllers are either implicit or explicit. The design of implicit dual controls, which we do not address in this paper, is based on the idea to approximate the Bellman equation of the problem. See \cite{mesbah_stochastic_2017} for a review on implicit dual control and \cite{bayard_implicit_2008} for an example of use of a particle filter in implicit dual control. On the contrary, in explicit dual control, one modifies the original problem to incorporate an explicit excitation of the system to maintain the dual effect. Explicit dual effect can be included in the optimization problem as a constraint or in the cost. Including it as a constraint generally leads to controllers with persistent excitation, constraints on information or approximation with scenario trees. Including it in the cost leads to integrated experiment design, where a measure of information is added to the original cost. In \cite{marafioti_persistently_2014}, a persistent excitation controller is presented but it is supposed that the system is linearizable which is not a suitable assumptions in our framework. Controllers with constraints on the information may lead to infeasibilty problems that are also hard to anticipate in nonlinear cases. In \cite{telen_study_2017}, this infeasibility issue is addressed in the deterministic framework. Scenario tree methods are computationally demanding and usually rely on Kalman filters inside the optimization problem, like in \cite{hanssen_scenario_2015}, \cite{subramanian_economic_2015} which are not adapted to multimodal cases. 

In integrated experiment design, a new term based on the Fisher information matrix (FIM) is added to the original cost as classically done in optimal design \cite{fedorov_model-oriented_2012}. However, if the original cost to minimize conditions the stability of the system, then the resulting trade off may destabilize it.

That is why, in this paper, we propose an output feedback dual SMPC based on integrated experiment design for information probing and on a Lyapunov constraint on the control to ensure stability. Output feedback is obtained by setting the initial condition of the optimization problem as particles from a particle filter that is able to deal with state estimation for arbitrary systems.

The paper is organized as follows. Section \ref{section_notation} recalls some important probabilistic notations. Section \ref{section_setup} recalls the bases of SMPC with imperfect information and describes our approach. Section \ref{section_par_fil} recalls the principle of particle filtering and section \ref{section_control_policy} describes the particular optimization problem that is solved in a receding horizon manner. In \cite{hokayem_stochastic_2012} and \cite{homer_output-feedback_2017}, stochastic output feedback  stability is proved but a separation principle was assumed. Actually, proving stochastic output feedback stability in the general case, when observability depends on the control, and when the only efficient estimators are particle filters, has never been proved to the best of our knowledge. However, we demonstrate the efficiency of our method on a challenging application in section \ref{application}.
\section{NOTATION}\label{section_notation}

Let $(\Omega, \mathcal{F}, P)$ be a probability space. In the following, random variables refer to $\mathcal{F}$-measurable functions defined on $\Omega$. For $i\in \mathbb{N}$, $\mathcal{B}(\mathbb{R}^i)$ denotes the set of Borel sets of $\mathbb{R}^i$. For a random variable X and a probability distribution, $X\sim p$ means that $p$ is the probability law of $X$. $P(\cdot\vert\cdot)$ and $E(\cdot\vert\cdot)$ denotes conditional probability and expectation. $P_x$ and $E_x$ denotes especially probability and expectation conditionally to $X_0=x$, for $x\in\mathbb{R}^i$. $P_p$ and $E_p$ denotes probability and expectation conditionally to $X_0\sim p$.
\section{PROBLEM SETUP} \label{section_setup}
\subsection{Description of the system}
We consider a discrete-time controlled stochastic dynamical system ${(X_k)}_{k \in \mathbb{N} }$ valued in $\mathbb{R}^{n_x}$ described by the following equation, $\forall k \in \mathbb{N}$:
\begin{align}
&X_{k+1}=f(X_{k},U_{k},{\xi}_{k}) \label{sys_dyn},
&X_0\sim p_0,
\end{align} where:
\begin{itemize}
{\setlength{\baselineskip}{1.2\baselineskip}
\item $p_0$ is the initial probability law.
\item 
${(U_k)}_{k \in \mathbb{N} }$ is the control process valued in $\mathbb{R}^{n_u}$.
\item ${(\xi_k)}_{k \in \mathbb{N}}$ are i.i.d. random variables valued in $\mathbb{R}^{n_{\xi}}$ distributed according to $p_\xi$. For each $ k \in \mathbb{N}$, $\xi_k$ represents an external disturbance on the dynamics.
\par}
\item f : $\mathbb{R}^{n_x}\times\mathbb{R}^{n_u}\times\mathbb{R}^{n_{\xi}} \longrightarrow \mathbb{R}^{n_x}$ is measurable.
\end{itemize}

We also assume that the state of the system is only available through some observations represented by a stochastic process ${(Y_k)}_{k \in \mathbb{N}}$ valued in $\mathbb{R}^{n_y}$ which verifies, $\forall k \in \mathbb{N} $:
\begin{align}
Y_k=h(X_k,\eta_k)\label{eq_obs},
\end{align}
where:
\begin{itemize}
\item ${(\eta_k)}_{k \in \mathbb{N} }$ are i.i.d. random variables valued in $\mathbb{R}^{n_{\eta}}$ distributed according to $p_\eta$. For each $ k \in \mathbb{N}$, $\eta_k$  represents an external disturbance on the observations. 
\item $h : \mathbb{R}^{n_x}\times\mathbb{R}^{n_{\eta}} \longrightarrow \mathbb{R}^{n_y}$ is measurable.
\end{itemize}
For $k \in \mathbb{N}$, we define the \textit{information vector} $I_k$ as follow:
\begin{align}
I_k&=(Y_0,U_0,\dots,Y_{k-1},U_{k-1},Y_k),\label{inf_vec}
\end{align}
From $I_k$, one can derive two important quantities in stochastic control with imperfect information that are the conditional distribution of $X_k$ given $I_k$ denoted by $\mu_k$ and the conditional distribution of $X_{k+i}$ given $(I_k,U_k,\dots,U_{k+i-1})$ for any $i \in \mathbb{N}^*$, denoted by $\mu_{k+i\vert k}$. Moreover, we denote by $K$ the Markov kernel defined by equation (\ref{sys_dyn}) and we assume that the conditional distribution defined by equation (\ref{eq_obs}) has a density with respect to the Lebesgue measure such that there exists a likelihood function denoted by $\rho$.
Therefore, for $k \in \mathbb{N}$, $A \in \mathcal{B}(\mathbb{R}^{n_x})$ and $B \in \mathcal{B}(\mathbb{R}^{n_y})$ :
\begin{align*}
    P(X_{k+1}\in A \vert X_k=x_k, U_k=u_k)&= K(A,x_k,u_k),\\
    P(Y_{k}\in B \vert X_k=x_k)&=\int_B \rho(y_k,x_k)\mathrm{d}y_k.
\end{align*}
$\mu_0$ is supposed to be known. Thus, $\forall (k,i) \in \mathbb{N}^2$,  $\mu_k$ and $\mu_{k+i\vert k}$ verify the following nonlinear filtering equations that can be summed up
to:
\begin{align}
    \mu_{k+1}&=F\left(\mu_k,Y_{k+1},U_k\right),\label{filtering_equation_red}\\
    \mu_{k+i+1\vert k}&=G\left(\mu_{k+i\vert k},U_{k+i}\right).\label{filtering_prop_red}
\end{align}
The expressions of $F$ and $G$ are classic and can be found in \cite{mesbah_stochastic_2017}. A control policy at time $k$ is then a map, denoted by $\pi_k$, that maps a conditional distribution $\mu_k$ to a control $U_k$. We denote a sequence of control policies by:
\begin{align*}
    \pi_{i:j}:=(\pi_i,\dots,\pi_j) ,\;\text{for}\; i\leq j.
\end{align*}

\subsection{Dual stochastic model predictive control}
Stochastic model predictive control is a widely used method for designing controllers in the presence of possibly unbounded disturbances in nonlinear dynamics as the one described in equation (\ref{sys_dyn}). It consists in solving a finite-horizon discrete stochastic optimal problem, to only apply the first control policy of the optimal sequence and to solve the problem again starting from the new state of the system. See \cite{mesbah_stochastic_2016} for a review on general stochastic MPC.

A MPC scheme is defined by the following features: 
\begin{itemize}
    \item A time horizon, $T\in \mathbb{N}^*$. 
    
    \item A family of set of constraints on the control, $\mathbb{U}_i\subset\mathbb{R}^{n_u}$, $\forall i = 0,..,T-1 $.
    \item $\forall i = 0,..,T-1 $, an instantaneous cost $g_i$: $\mathbb{R}^{n_x}\times\mathbb{R}^{n_u}\times\mathbb{R}^{n_{\xi}} \longrightarrow \mathbb{R}$  and a final cost $g_T$: $\mathbb{R}^{n_x} \longrightarrow \mathbb{R}$.
     \end{itemize}
     
     The choice of the objective functions $g_i$ and $g_T$ and of the control constraints $\mathbb{U}_i$ is a matter of modelling. Indeed, one often has an economic cost to minimize that comes from practical considerations, such as a price or a fuel consumption. Moreover, the control must be designed to attain some target in the state space: it is the guiding problem. There are two classical ways to address this issue in MPC:
     \begin{itemize}
         \item Adding a new term to the cost to enforce stability in some sens. In this case, the general cost is decomposed in the following way, for $i = 0,..,T$:
    \begin{align}
        g_i=g_{i}^{\textrm{stab}}+g_{i}^{\textrm{eco}},\label{inst_cost}
    \end{align}
     Thus, $g_i$ and $g_T$ realize a compromise between convergence and economic costs. For instance, in the LQG case, $g_i(x,u,\xi)=x^T M_xx+u^TM_u u$ where $M_x$ and $M_u$ are positive definite matrices. The first term drives the state of the system to zero and the second term penalizes high controls. The compromise is dealt with by tuning the matrices $M_x$ and $M_u$.
     \item Adding a drift constraint on the first control, $U_0$, that enforces the decreasing of some Lyapunov-like function, during the first time step only, such that, for $i = 0,..,T$:
     \begin{align}
        &g_i=g_{i}^{\textrm{eco}},\label{inst_cost_eco}\\
        &\text{a negative drift condition on } U_0. \;\notag
    \end{align}
     Actually, since only $U_0$ is applied on the system, the Lyapunov function decreases along the whole trajectory and then stability is obtained. It is also known as Lyapunov Economic MPC, see \cite{ellis_tutorial_2014} for a review in the deterministic setting. In the stochastic setting, it has been applied with output feedback for continuous-time nonlinear systems in \cite{homer_output-feedback_2017} and for a discrete-time linear system with bounded controls in \cite{hokayem_stochastic_2012}. 
     \end{itemize}

In the presence of partial information represented by equation (\ref{eq_obs}), the theoretical receding horizon problem to solve is more complex than in the perfect information case, as the the new state of the system is $\mu_k$ \cite{bertsekas_dynamic_2005}, and it evolves according to equation (\ref{filtering_equation_red}). Thus, the problem to solve at each time k denoted by $(P^k_C)$ can be written as follows:

\begin{equation*}
\begin{array}{rrclcc}
\displaystyle \underset{\pi_{0:T-1}}{\text{min}} & \multicolumn{3}{l}{E_{\mu_k}\left[\sum_{i=0}^{T-1}E[ g_i(X_i,U_i,\xi_i)\vert I_i] + E[g_T(X_T)\vert I_T]\right]} \\
\textrm{s.t.} & \tilde{\mu}_{i+1} & = & F\left(\tilde{\mu}_{i},Y_{i+1},U_i\right), \\
& U_i & = & \pi_i(\tilde{\mu}_{i}), \; \\
& U_i & \in & \mathbb{U}_i, \; \forall i = 0,..,T-1 , \\
& \tilde{\mu}_{0} & = & \mu_k.
\end{array}
\end{equation*}

In $(P^k_C)$, the conditional expectation w.r.t. the available information $I_i$  emphasizes the fact that the observations, are propagated from $k$ to $k+T$ and depend on the control. Therefore, as recalled in \cite{mesbah_stochastic_2017} and \cite{flayac_nonlinear_2017}, one can see from the Dynamic Programming Principle that the optimal control policies of $(P^k_C)$ generally have the dual effect property \cite{bar-shalom_dual_1974}. As optimal policies are intractable in the general nonlinear case due to the curse of dimensionality, suboptimal policies that are constructed instead. They are called dual controllers as they preserve the dual effect property by: 
\begin{itemize}
    \item controlling the system in the standard way defined by equation (\ref{inst_cost}) or (\ref{inst_cost_eco}).
    \item explicitely or implicitely probing information to improve the quality of the future observations.
\end{itemize}
 In this article, we focus on explicit dual controllers that take information into account inside the cost. This technique is called integrated experiment design. In fact, the new cost denoted by $g_i^{\textrm{ex}}$ realizes a trade off between the original costs $g_i$ and $g_T$ and a measure of information denoted by $g_i^{\textrm{info}}$ such that, for $i = 0,..,T$:
\begin{align}
    g_i^{\textrm{ex}}=g_{i}+g_{i}^{\textrm{info}},\label{inst_cost_inf}
\end{align}

An important issue appears when, additionally, the guiding objective is incorporated in the cost as described in equation ($\ref{inst_cost}$). Consequently, (\ref{inst_cost_inf}) represents the sum of three terms that usually contain weights that must be tuned. By construction, the weights affect the convergence of the system to the guiding goal which makes their tuning difficult. Indeed, if too much importance is given to $g_{i}^{\textrm{info}}$ then probing will be efficient but the system may get stuck far from the target. On the contrary, if too much importance is given to $g_{i}^{\textrm{stab}}$ then the probing effect will not be sufficient and output feedback performance may be poor. It is difficult to know which case will occur a priori depending of the value of the weights because the optimal costs are not known at first sight.

 That is why, in this paper, we present an output-feedback explicit dual controller based on the minimization of the Fisher information matrix subject to a negative mean drift condition for nonlinear discrete-time systems coupled with a particle filter for state estimation. In other words, we have chosen to model our problem with equations (\ref{inst_cost_eco}) and (\ref{inst_cost_inf}).   
 
 Furthermore, in explicit dual control, one does not have to propagate the information inside the dynamics with equation (\ref{filtering_equation_red}) anymore because it is dealt with in a different way. It is sufficient to propagate $\mu_{i\vert0}$ with equation (\ref{filtering_prop_red}) which is much simpler. Besides, one wants to solve in practice a finite dimensional optimization problem so one looks for control values and not policies anymore. Then, the general optimization problem associated with our method, denoted by $(P^k_{\textrm{ex}})$, can be written as follows: 

\begin{equation*}
\begin{array}{rrclcc}
\displaystyle \underset{u_{0:T-1}}{\text{min}} & \multicolumn{3}{l}{E_{\mu_k}\left[E[\sum_{i=0}^{T-1} g_i^{\textrm{ex}}(X_i,u_i,\xi_i) + g_T^{\textrm{ex}}(X_T)\vert I_0]\right]} \\
\textrm{s.t.} &  \tilde{\mu}_{i+1\vert 0}&=&G\left(\tilde{\mu}_{i\vert 0},u_{i}\right), \\
& u_i & \in & \mathbb{U}_i, \; \forall i = 0,..,T-1 , \\
& \tilde{\mu}_{0} & = & \mu_k.\\
&\text{negative}&\text{drift}&\text{condition on}\;u_0.\\
\end{array}
\end{equation*}

\section{STATE ESTIMATION}\label{section_par_fil}

From the definition of $(P^k_{\textrm{ex}})$, it is clear that dual SMPC requires the computation of $\mu_k$. However, in the general nonlinear case, $\mu_k$ cannot be computed explicitly so approximation methods are needed. Kalman filters are widespread and easy-to-compute approximations of the posterior distribution but they may fail in the presence of high nonlinearities and multimodality. That is why we use particle filters that are known to handle these difficulties at the price of a higher computational cost. A particle filter approximates the posterior distribution $\mu_k$ by a set of N particles, ${\left(x^{l}_{k}\right)}_{l= 1,..,N }$ valued in $\mathbb{R}^{n_x}$, associated with nonnegative and normalized weights ${\left({\omega}^{l}_{k}\right)}_{l=1,..,N }$. The approximate distribution, denoted by, $\mu_k^N$, is then defined as follows:
\begin{align}
      \mu_k^N=\textstyle\sum_{l=1}^{N}\omega_k^{l}\delta_{x^{l}_k} \label{particle_filter_dirac},
\end{align}

Algorithm \ref{particle_filter} described the steps of a particle filtering algorithm  with adaptive resampling. Furthermore, we define the conditional expectation of the state, denoted by $\hat{x}_k$ w.r.t. the information $I_k$ and its particle approximation, the empirical mean of the filter $\hat{x}_k^N$ i.e:
\begin{align*}
    &\hat{x}_k=E[X_k\vert I_k],
     &\hat{x}_k^N=\frac{1}{N} \textstyle\sum_{l=1}^N \omega^{l}_{k} x^{l}_{k}.
\end{align*}

\section{CONTROL POLICY}\label{section_control_policy}
The goal of this section is to present our new explicit dual receding horizon control scheme. The main idea is to prioritize the guiding goal by adding a stabilizing constraint on the first control of the finite horizon optimization problem. This constraint consists in forcing the decreasing of a Lyapunov-like function during the first time step starting from the current estimator the state. It is equivalent to stabilizing an estimator of the state. The success of such a technique depends highly on the estimation error. Indeed, if the estimation error is high or even diverges, driving the estimator to the target does not imply driving the true state to the target. It is a well known problem  in the deterministic setting \cite{andrieu_unifying_2009}. In \cite{homer_output-feedback_2017} and \cite{hokayem_stochastic_2012}, output feedback stability is proved but it was assumed that the estimation error converges uniformly w.r.t. the control. This assumption is not realistic anymore in our setting because the observations do depend on the control and the estimation error may diverge if the dual effect is not taken into account. Accordingly, our method consists in solving a version of  $(P^k_{\textrm{ex}})$ where $g_i^{\textrm{ex}}$ is based on the Fisher information matrix and the drift condition comes from the theory of Markov chains stability.

%

\subsection{Foster-Lyapunov drift in case of perfect information}
Perfect information is met when full knowledge of the state $X_k$ is available. In the time homogeneous case, the control is computed from a state feedback control policy i.e. a measurable maps, denoted by $\alpha$, that maps a state $X_k$ to a control $U_k$. Thus, for any $\alpha$, one can define the corresponding closed loop system as follows, $\forall k \in \mathbb{N}$:
\begin{align}
&X_{k+1}=f(X_{k},\alpha(X_k),{\xi}_{k}) \label{sys_dyn_perfect_info},
&X_0\sim p_0.
\end{align}
 A state feedback control policy is said to be \textit{admissible} if $\forall x\in \mathbb{R}^n, \alpha(x)\in \mathbb{U}_0$.
Therefore, equation (\ref{sys_dyn_perfect_info}) defines also a time-homogeneous Markov chain whose stability can be studied via the classical theory of negative drifts conditions discussed in \cite{meyn_markov_2009} and recalled in \cite{chatterjee_stability_2015}. In proposition \ref{prop_drift}, we focus on \textit{geometric drifts conditions} that are closely related to Lyapunov conditions for exponential stability for continuous-time processes.
\begin{algorithm}[htbp]
{\setlength{\baselineskip}{1\baselineskip}
\caption{Particle filter with adaptive resampling}\label{particle_filter}
\begin{algorithmic}[1]
\State Create a sample of $N$ particles $x^{l}_{0}$ according to the law $\mu_0$ and initialize the weights $\omega_0^{l}$ with $\frac{1}{N}$
\For{$i=0,1,2\dots$} 
\State \textbf{Prediction:}
\State Given a control $u_i$ and the particles ${\left(x^{l}_{i}\right)}_{l= 1,..,N }$, compute the predicted particles by drawing samples from $K$ i.e.
\begin{align*}
    x^{l}_{i+1\vert i}\sim K(\mathrm{d}x_{i+1\vert i},x_i^l,u_i), \;\text{for}\; l=1,..,N.
\end{align*}
\State \textbf{Correction}:
\State Get the new observation $Y_{i+1}$
\State Compute the updated weights ${\left(\tilde{\omega}^{l}_{i}\right)}_{l=1,..,N }$ thanks to the likelihood function $\rho$:
\begin{align*}
    \tilde{\omega}^{l}_{i}={\omega}^{l}_{i}\rho(Y_{i+1},x^{l}_{i+1\vert i})
\end{align*}
\If{\textit{Resampling}}
\State Draw the \textit{a posteriori} particles ${\left(x^{l}_{i+1}\right)}_{l= 1,..,N }$ from the set ${\left(x^{l}_{i+1\vert i}\right)}_{l= 1,..,N }$ and ${\left(\tilde{\omega}^{l}_{i}\right)}_{l=1,..,N }$ using a resampling technique and set $\omega^{l}_{i+1}=\frac{1}{N}$
\Else 
\State Set $x^{l}_{i+1}=x^{l}_{i+1\vert i}$ and ${\omega}^{l}_{i+1}=\tilde{\omega}^{l}_{i}/\sum_{l=1}^N\tilde{\omega}^{l}_{i} $
\EndIf
\EndFor
\end{algorithmic}
\par}
\end{algorithm}
\newtheorem{prop_drift}{Proposition}

\begin{prop_drift} \label{prop_drift}
Suppose that there exist $b>0$ and $\lambda_{min}\in [0,1[$, a measurable function $V:\mathbb{R}^{n_x} \longrightarrow [0,+\infty[$, a compact set $C\subset\mathbb{R}^{n_x} $ and an admissible state feedback control policy $\alpha$ such that $E[V(f(x,\alpha(x),{\xi}_{0})]\leq \lambda_{min} V(x)$, $\forall x \notin C $ and $\textrm{sup}_{x \in C}E_x[V(X_1)]=b$

Then, $\forall \lambda\in [\lambda_{min},1[$, 
\begin{align*}
    E_x[V(X_k)]\leq \lambda^k V(x)+b(1-\lambda)^{-1}, \forall k\in \mathbb{N},\;\forall x\in \mathbb{R}^n, 
\end{align*}
 where $X_k$ is computed with equation (\ref{sys_dyn_perfect_info}). 
\end{prop_drift}

Proposition \ref{prop_drift} is a slightly different reformulation of proposition 1 in \cite{chatterjee_stability_2015}, but its proof follows from the one in \cite{chatterjee_stability_2015}. In particular, we consider that the parameter $\lambda$ can be chosen arbitrarily in $[\lambda_{min},1[$. In practice, $\lambda$ is a parameter to tune that determines the convergence speed of the system. Moreover, as explained in \cite{chatterjee_stability_2015}, if proposition \ref{prop_drift} is verified for a norm-like function V then for $r>0$, $P_x(\Vert X_k \Vert>r)$ decreases as the inverse of V so the distribution of state concentrates itself around $0$. That is why, in the rest of the paper, we suppose that our guiding goal is to drive the system (\ref{sys_dyn}) to 0. To do that, we also suppose that the assumptions in proposition \ref{prop_drift} are fulfilled and notably that the system (\ref{sys_dyn}) can be stabilized with perfect information with some admissible state feedback control policy $\alpha$.

\subsection{Receding horizon policy}
In our dual MPC scheme, the information is quantified by the Fisher information Matrix (FIM), denoted by $J_i$. Recursive computation of the FIM can be found in \cite{tichavsky_posterior_1998}. Thus, explicit dual effect is created by minimizing some functions of the FIM, $g_{i}^{\textrm{fish}}$ and $g_{T}^{\textrm{fish}}$ such that, in $(P^k_{\textrm{ex}})$ :
\begin{align}
     g_i^{\textrm{ex}}(X_i,U_i,\xi_i)&=g_{i}^{\textrm{eco}}(X_i,U_i,\xi_i)+g_{i}^{\textrm{fish}}(J_i),\label{inst_cost_exp}
\end{align}
In $(P^k_{\textrm{ex}})$, the drift condition is taken from proposition \ref{prop_drift} and, applied to $\hat{x}_k$ only when $\hat{x}_k \notin C$ such that, for some $ \lambda\in [\lambda_{min},1[$:
\begin{align*}
    E_{\hat{x}_k}[V(f(\hat{x}_k,u_0,{\xi}_{0})]\leq \lambda V(\hat{x}_k),\;\textrm{when}\; \hat{x}_k \notin C.
\end{align*}

It is important to notice that the admissibility of the drift constraint is guaranteed by the existence of the stabilizing admissible state feedback policy $\alpha$.

A good approximation technique to approach $({P}_{\textrm{ex}}^k)$ in practice is the scenario approach \cite{calafiore_scenario_2006}. It appears that, in output feedback MPC, there is a synergy between the scenario approach and particle filtering. In fact, the initial condition for each independent scenario is chosen as a particle from the current set of particles. It improves global performance compared to a similar technique involving a Kalman filer in which, the initials conditions are always drawn according to a unimodal law. This method was already used in \cite{sehr_particle_2016} and \cite{flayac_nonlinear_2017}. The approximation, denoted by $\left({P}_{\textrm{ex}}^{k,N}\right)$, can be defined, $\forall \lambda\in [\lambda_{min},1[$, by:
\begin{equation*}
\begin{array}{rrclcc}
\displaystyle \underset{u_{0}\cdots u_{T-1}}{\text{min}}\hspace{-0.3cm} & \multicolumn{3}{l}{\sum_{l=1}^{N_s}\omega_k^{l}\left(\sum_{i=0}^{T-1} {g}_i^{\textrm{ex}}\left(X^{l}_{i},u_i,\xi_i^l \right) + {g}_T^{\textrm{ex}}\left(X^{l}_{T}\right)\right)} \\
\\
 \textrm{s.t.}&X^{l}_{i+1} & = &f(X^{l}_{i},u_{i},{\xi}^{l}_{i}),\\
 &u_i & \in & \mathbb{U}_i,\\
 & X^{l}_{0} & = & {x}^{l}_{k},\;\forall l =1,..,N_s,\; \forall i =  0,..,T-1 , \\
\end{array}
\end{equation*}
\begin{equation*}
    \frac{1}{N_{{dr}}}\textstyle\sum_{\ell=1}^{N_{{dr}}}V(f(\hat{x}_k^N,u_0,\tilde{\xi}_{0}^{\ell})\leq \lambda V(\hat{x}_k^N),\;\textrm{when}\; \hat{x}_k^N \notin C,\\
\end{equation*}
where:
\begin{itemize}
    \item $N_s<N$ is the number of scenarios considered. It is supposed to be less than the number of particles for computational reasons, so $N_s$ particles must be extracted from the original set. 
    \item $N_{dr}\in \mathbb{N}^*$ is the size of the sample used to approximate the expectation in the drift constraint.
    \item ${({\xi}^{l}_{i})}_{ i =  0,..,T-1, l =  1,..,N_s}$ and ${(\tilde{\xi}_{0}^{\ell})}_{\ell=1,..,N_{dr}}$ are i.i.d. random variables sampled from $p_{\xi}$.
    \end{itemize}

Additionally, when $\hat{x}_k^N \in C$ the drift condition is not necessarily feasible so we decided to apply the stabilizing policy $\alpha$ to $\hat{x}_k^N$, such that $U_k=\alpha(\hat{x}_k^N)$, as it is done for classical certainty equivalence controllers. This means that if $\hat{x}_k^N$ enters $C$ there is no probing anymore. It is not a problem as one can see $C$ as the target set.
\subsection{Global algorithm}
The complete Output feedback algorithm is summarized in Algorithm \ref{algo_fisher}. This method has two main advantages compared to the one we presented in \cite{flayac_nonlinear_2017}.

 First, in \cite{flayac_nonlinear_2017}, both the stability properties and the information were incorporated as two terms of the cost. Therefore, the weights between the terms were difficult to tune and, especially in receding horizon control where they deeply influence the convergence of the system. In particular, it was compulsory to decrease the weights on $g_{i}^{\textrm{fish}}$ and $g_{T}^{\textrm{fish}}$ with time otherwise the system converged to a point that was far from the target. Moreover, tuning the decreasing of this weight was also complicated a priori. In our new method, the stability properties are much less influenced by the cost because of the drift condition which is a constraint. The most important parameter to tune  is $\lambda$, and, in principle, it influences only the convergence speed of the system and not its qualitative properties of stability.
 
  Secondly, from a numerical point of view, as stated in \cite{homer_output-feedback_2017}, the stability properties are contained in one constraint and are easier to achieve in practice. Indeed, the convergence of the system depends little on the quality of the solution of the optimization problem and much more on its admissibility which is easy to obtain with classical solvers.

\begin{algorithm}[h]
{\setlength{\baselineskip}{1\baselineskip}
\caption{Fisher/Lyapunov Output Feedback Control}\label{algo_fisher}
\begin{algorithmic}[1]
\State Create a sample of $N$ particles ${x}^{l}_{0}$ according to the law $\mu_0$ and initialize the weights $\omega_0^{l}$ with $\frac{1}{N}$.
\For{$k=0,1,2,\dots$}
\If{$\hat{x}_k^N \notin C$}
\State Solve $\left({P}_{\textrm{ex}}^{k,N}\right)$ starting from the particles ${x}^{l}_{k}$ and the weights $\omega_k^{l}$.
\State Get an optimal sequence $(u_0^*,\dots,u_{T-1}^*)$.
\State Set $U_k=u_0^*$.
\Else 
\State Set $U_k=\alpha(\hat{x}_k^N)$
\EndIf
\State Compute the \textit{a posteriori} particles $x^{l}_{k+1}$ and weights $\omega_{k+1}^{l}$ given $U_k$ according to Algorithm \ref{particle_filter}
\EndFor
\end{algorithmic}
\par}
\end{algorithm}
\section{APPLICATION AND NUMERICAL RESULTS}\label{application}
\subsection{Description of the application}
Algorithm \ref{algo_fisher} has been applied to the guidance and localization of a drone by terrain-based navigation (TBN). The guiding goal is to drive a drone in a 3D space from an uncertain initial condition $X_0$ to a compact set centered around $0$. If the original target is not 0 then a translation can be made to center the problem around 0. It is assumed that the drone can be described by 3 positions and 3 speeds, $X_k={(x_k,y_k,z_k,v^x_k,v^y_k,v^z_k)}^T$ and the control of 3 accelerations $U_k={(u^x_k,u^y_k,u^z_k)}^T$. $(x_k,y_k)$ represents a horizontal position and $z_k$ an altitude. Its dynamics is linear with bounded controls such that, for $k\in \mathbb{N}$ :
\begin{align}
        &X_{k+1}=AX_k+BU_k+\xi_k, &\Vert U_k \Vert\leq U_{max}, \label{sys_lin}
\end{align}
where $U_{max}>0$ and $A\in \mathbb{R}^{n_x\times n_x}$ and $B\in \mathbb{R}^{n_x\times n_u}$ correspond to a discrete-time second order system with damping on the speed. Because of the constraint on the control, classical linear controllers are not feasible and the closed loop system must be nonlinear.
The guiding problem is addressed with a drift constraint taken from proposition 8 in \cite{chatterjee_stability_2015} such that $V(X)=e^{\Vert X \Vert}$ and C is a ball centered around 0 whose radius depends only on the disturbances of the dynamics. To guarantee the existence of an admissible stabilizing state feedback policy, $\alpha$, the maximum control $U_{max}$ must be sufficiently high. See \cite{chatterjee_stability_2015} for the precise definition of $\alpha$, $U_{max}$ and $C$. 

 The speed is supposed to be measured. Concerning the position, the paradigm of TBN is that only a measure of the difference between the altitude of the drone, $z_k$, and the altitude of the corresponding vertical point on the ground is available. The ground is represented by a map, $h_{m}$, that maps a horizontal position $(x,y)$ to the corresponding height of the terrain. In practice, $h_{m}$ is determined by a smooth interpolation of data points so it is highly nonlinear. Therefore, the observation equation is: 
 \begin{align*}
      h(X_k,\eta_k)={(z_k-h_{m}(x_k,y_k),v^x_k,v^y_k,v^z_k)}^T+\eta_k. 
 \end{align*}
 
 Because of the potential multimodality of $h_m$, state estimation must be dealt with by a particle filter. Moreover, it appears very naturally that dual control is required in this application. Indeed, the quality of the observations depends on the control and more precisely on the area that is flied over by the drone. Let us assume that the drone flies over a flat area with constant altitude then one measurement of height matches a whole horizontal area and the estimation error on $(x_k,y_k)$ is of the order of magnitude of the size of the area. On the contrary, if the drone flies over a rough terrain, then one measurement of height corresponds to a smaller area on the ground and the estimation error is reduced. Finally, an intuitive good dual control strategy is to go toward the target avoiding to fly over uninformative areas.
\subsection{Numerical results} 
\begin{figure}[h]
\begin{center}
\includegraphics[scale=0.67]{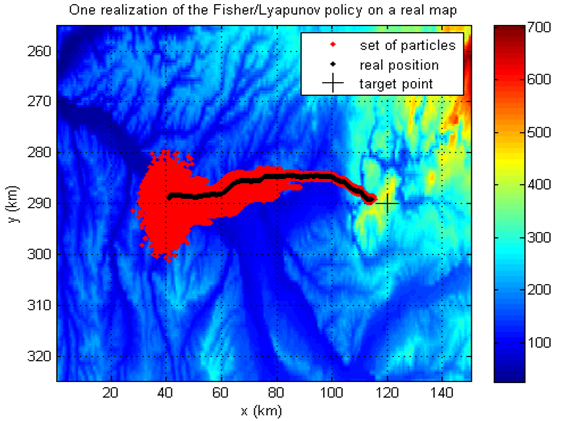} 
\end{center}
\caption{Plot of one trajectory obtained by Fisher/Lyapunov particle control and of the particles from the particle filter}
\label{fisher_traj}
\end{figure}

\begin{figure}[h]
\begin{center}
\includegraphics[height=6.5cm,width=0.5\textwidth]{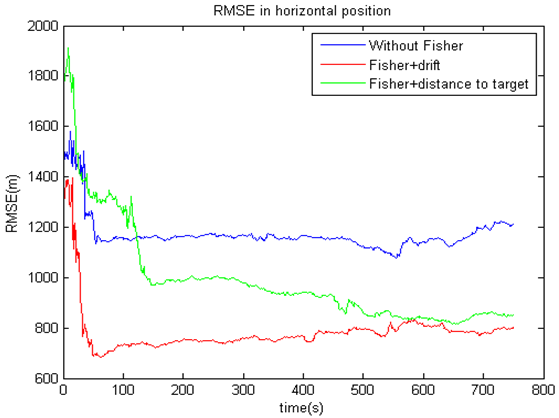}
\end{center}
\caption{Plot of the RMSE in horizontal position for 3 policies: without the FIM (blue), with the FIM and the drift (red), with the FIM and the distance to the target (green)}
\label{RMSE}
\end{figure}

Figure \ref{fisher_traj} represents the horizontal projection of a trajectory obtained by algorithm \ref{algo_fisher} with the terrain map in its background. The scenario has been chosen such that, if the system goes in straight line to the target then it flies over a flat area. We can see that in figure \ref{fisher_traj} that, as expected, the system makes a detour to avoid the flat area, so that the particles tighten around the true position, and finally reaches the target. In this run, the algorithm has been stopped when the target enters the compact C.
Figure \ref{RMSE} shows RMSE in horizontal position after 30 Monte Carlo simulations of three different policies: a constrained LQ MPC with no information probing, the policy described in this article and the policy described in \cite{flayac_nonlinear_2017} where stability was forced by the minimization of the distance to the target. 
%


The simulations were run in MATLAB and the optimization problems were solved using the modelling language AMPL and the solver Ipopt.
%
\section{CONCLUSION}
In this paper, we have presented a new explicit dual output feedback stochastic MPC for nonlinear systems. Its principle is to choose at each time step the controls that maximizes the information over the control that forces the mean decreasing of some Lyapunov-like function for discrete-time nonlinear stochastic systems. In particular, it does not involve penalization of the guiding goal in the cost which is a classical feature of MPC. Output feedback is obtained by coupling the resolution of an optimization problem with a particle filter. The method is applied to terrain aided navigation and appears to be easier to tune that one of our previous method. 
 


\bibliographystyle{unsrt}
\bibliography{articles_these}

\end{document}